%% file: document.tex
\documentclass[a4paper]{amsart}
\usepackage{amsmath, amsthm, amssymb, amsbsy}

\usepackage[pdftex]{graphicx}
\usepackage{tikz}
\usepackage{enumerate}
\usepackage{listings}
\usepackage{caption}
\DeclareCaptionFont{white}{\color{white}}
\DeclareCaptionFormat{listing}{\colorbox{gray}{\parbox{\textwidth}{#1#2#3}}}
\usepackage{pgflibraryshapes}
\usetikzlibrary{shapes,arrows,decorations,backgrounds,automata}
\usetikzlibrary{mindmap,trees}
\usetikzlibrary{decorations.pathreplacing}
\usetikzlibrary{plotmarks,patterns}
\usetikzlibrary{3d,calc}

\graphicspath{{./}}

\begin{document}

\title{A Homotopy Method for Large-Scale Multi-Objective Optimization}

\author[Adelmann]{Andreas Adelmann}
\email[Adelmann]{andreas.adelmann@psi.ch}
\author[Arbenz]{Peter Arbenz}
\author[Foster]{Andrew Foster}
\author[Ineichen]{Yves Ineichen}


\input{Abstract.tex}

\maketitle

\section{Introduction}
\input{Introduction.tex}

\section{Motivating Example}
\label{sec:Motivation}
\input{Motivation.tex}

\section{Multi-Dimensional Generalization}
\label{sec:Generalization}
\input{Generalization.tex}

\section{Results}
\input{Results.tex}


\section{Conclusions}
\input{Conclusions.tex}

\input{Bibliography.tex}

\end{document}

%% file: Abstract.tex
\begin{abstract}

A homotopy method for multi-objective optimization that produces uniformly
sampled Pareto fronts by construction is presented. While the algorithm is general, of
particular interest is application to simulation-based engineering optimization
problems where economy of function evaluations, smoothness of result, and
time-to-solution are critical. The presented algorithm achieves an order of
magnitude improvement over other geometrically motivated methods, like Normal
Boundary Intersection and Normal Constraint, with respect to solution evenness
for similar computational expense. Furthermore, the resulting uniformity of
solutions extends even to more difficult problems, such as those appearing in
common Evolutionary Algorithm test cases. 


\end{abstract}

%% file: Introduction.tex
The problem of scalar function minimization is ubiquitous in modern engineering
sciences and a number of algorithms exist that exploit various functional
characteristics to reach timely solutions. Far more difficult is the problem of
vector function minimization, where a of number different, often conflicting,
objectives are optimized in such a way as to strike a balance that pleases the
end-user. The presence of multiple output dimensions requires a generalization
of our concept of optimality, demands significantly more computational effort,
as well as compounds the inherent difficulty of the problem.

\subsection{Problem Statement}
Precisely stated, and adopting the notation of \cite{Marler2004}, we seek to
address the problem
\begin{eqnarray}
\label{prob:MOOP}
\min_{\textbf{x} \in D} \textbf{F}(\textbf{x}) &=& [F_1(\textbf{x}),
F_2(\textbf{x}), ..., F_k(\textbf{x})]^T \\
\text{subject to}  & & \nonumber \\
g_j(\textbf{x}) & \le 0 & \forall j \in [1, m] \nonumber \\
h_l(\textbf{x}) &=  0 & \forall l \in [1, e] \nonumber
\end{eqnarray}
where the objective functions, $F_i$, are bounded and defined over the set
$D \in \mathbb{R}^n$.  Here, $\textbf{x}$ is a vector of independent design
variables of length $n$ and contained in $D$.  In this formulation, the set of
objective functions, $ \left\{ F_i(\textbf{x}) \colon \forall i \in [1, k]
\right\} $, form a map
\begin{equation}
\textbf{F} \colon \mathbb{R}^n \to \mathbb{R}^k
\end{equation}
from design space to objective space. The set of all $\textbf{x} \in D$
satisfying the inequality constraints, $g_j(\textbf{x}) \le 0$, and the equality
constraints, $h_l(\textbf{x}) = 0$,
\begin{equation}
X \colon= \left\{ \textbf{x} \in D \colon  \begin{matrix} 
g_j(\textbf{x}) \le 0 & \forall j \in [1, m] \\
h_l(\textbf{x}) = 0 & \forall l \in [1, e]
\end{matrix} \right\}
\end{equation}
is called the feasible set.  It's image in objective space
\begin{equation}
Z \colon= \left\{ \textbf{F}(\textbf{x}) \colon \textbf{x} \in X \right\}
\end{equation}
is denoted the attainable set.

In the general case, the $\textbf{F}_i$ functions can not be simultaneously
optimized and compromise solutions must be considered.  Therefore, our goal is
to discover locally Pareto optimal points, which are solutions that can not be
improved in all objectives simultaneously when compared to its neighbors in the
feasible set.  The set of all Pareto optimal points (when viewed in objective
space) is called the Pareto front and explicitly maps the optimal values that
can be attained when considering conflicting objectives.

In this paper we propose a method for computing high-quality discrete
approximations of Pareto fronts.
We begin by reviewing existing methods and commenting on their ammenability to
large problems with many objectives.  We then review previous work on a
particular homotopy method and demonstrate an alternative formulation.  After
explicitly generalizing this approach to problems arbitary dimension, we
compare our formulation to methods using three test problems from the
literature.  Finally, we present an application of this method to a problem in
the field of particle accelerator design and discuss performance related to various
benchmarks.  We conclude by mentioning futher improvements and future research
directions.

\section{Related Work}

As a result of the ubiquity of problems that fit the formulation given in
\eqref{prob:MOOP}, a variety of methods have been devised for sampling
high-dimensional Pareto fronts.

\subsection{Scalarization}

One simple method involves explicitly imposing an order on the objective space. 
This can be done, for example, by introducing a suitable scalarization function
projecting the objective space onto, and using the well-ordered property of,
the real numbers to define a solution. This effectively reduces the task of
\eqref{prob:MOOP} to a scalar minimization problem amenable to solution via an
number of established numerical optimization techniques.
\begin{equation}
\label{prob:SCALAR}
\begin{matrix} 
\min_{\textbf{x} \in X} s( \textbf{F}(\textbf{x}) ; \boldsymbol\lambda) 
 \mbox{ 
with } 
s(\textbf{F} ; \boldsymbol\lambda ) \colon \mathbb{R}^k \to \mathbb{R}
\end{matrix}
\end{equation}
While the approach is straightforward, the mapping between the scalarization
parameters ($\boldsymbol\lambda$) and corresponding solution's location on the
Pareto front generally is not, and fronts generated by blindly probing the 
$\boldsymbol\lambda$-space often suffer from highly irregular sampling
\cite{Marler2004}.

\subsection{Evolutionary Algorithms}

Rather than attempting to break the overall sampling task into a series of
subproblems, as in the proposed scalarization approach, population-based
methods, generally speaking, attempt to solve the problem as a whole.  This is
done, typically, by generating sets of mutually non-dominated
solutions using Evolutionary Algorithms (EAs).  A number of specific heuristics
have been proposed (see \cite{Zitzler1999}, \cite{Marler2004}, \cite{Deb2002}, and \cite{Nebro2009})
using a variety of techniques to balance the discovery of non-dominated solutions with diversity
criteria to discourage solution crowding.

While EAs are capable of performing well on even the most difficult problems,
some key aspects suggest that they are not ideal for all objective types. First,
as in \cite{Deb2002}, solutions are evaluated based on their relative dominance
to other solutions in the population.  Unlike scalarization methods that provide
some necessary or sufficient conditions for optimality of solutions, this does
not actually guarantee Pareto optimality of the resulting solution candidates.
Additionally, the nature of the non-dominating fitness criteria inherently
limits the scalability of this class of algorithms as it involves comparing
elements of the population against each other, naturally implying some
super-linear running time behavior with increasing population size
\cite{Deb2002, Jensen2003}. While parallelizing the evaluation step
\cite{Durillo2008} or distributing various objective space regions to different
processors \cite{Deb2003} can partially mitigate this effect, population size
must still increase exponentially with objective space dimension to maintain
an evenly sampled Pareto front.  Furthermore, increasing the number of objectives
distorts the utility of the dominance ranking criterion as the hypersurface-area
to hypervolume effect permits a higher proportion of non-dominated solutions and
weakens selection pressure \cite{Ishibuchi2008}.

\subsection{Geometric Methods}

Rather than rely on a flood of function evaluations to explore the solution
space, a final class of methods uses the scalarization approach, but exploits
geometrical arguments to further restrict the attainable set for each
optimization subproblem.  Practically, this is accomplished by adding auxiliary
constraint functions to steer the scalar optimization subroutines towards
desirable solutions and deliver an evenly sampled front.  The Normal Boundary
Intersection (NBI) \cite{Das1998} and the Normal Constraint (NC) \cite{Messac2004} method
essentially restrict each optimization subproblem to consider only solutions
that lie on the normal vector emanating from the convex hull of the individual
minimizers (CHIM).

\begin{figure}[h!]
  \centering
  \input{./CHIM.tex}
  
  \caption{In the NBI and NC methods, evenly spaced points on the CHIM are
  projected to the Pareto front.  Unfortunately, as a result of the front
  geometry, the resulting sample points are not always evenly spaced.  The
  fact that $d_1 \ne d_2$ illustrates this effect.}
  \label{fig:CHIM}
\end{figure}
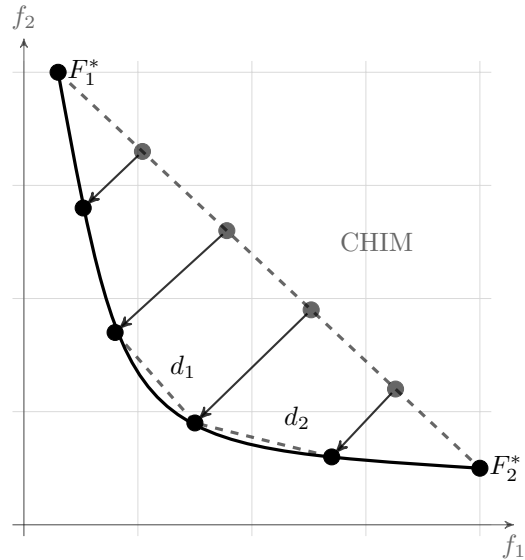

These methods make intuitive sense, as well as produce well sampled fronts. 
However, by restricting the attainable set to the CHIM normal,
these approaches are only capable of capturing Pareto points that lie within
the projection of this hull.  Furthermore, these points are only equally
distributed when viewed from the CHIM frame, as seen in figure \ref{fig:CHIM}.

The specific algorithms also exhibit some scalability bottlenecks, especially
when increasing the number of objectives, $k$.  Since the boundaries of the
Pareto front are not know a priori, it is often difficult to decide where to
place solution points on the CHIM.  In \cite{S.Motta2012} and
\cite{Mueller-Gritschneder2009}, it is suggested to divide the solution
algorithm into two phases; the first samples the trade-offs between all
combinations of individual objectives to establish the boundaries of the
complete front.  The second phase involves solving a set of scalar optimization
subproblems to find Pareto optimal points on the interior of this boundary.  The
first phase, consists of $\sum_{i=2}^{k-1} \binom{k-1}{i}$ ``sub-fronts'' to
sample before the domain is properly bounded and sampling on the interior can
begin.  Furthermore, the recursive nature of sub-front sampling, such as this, 
limits the scalability to only $i$ parallel sub-front's per step of the initial
phase.

%% file: CHIM.tex
  \begin{tikzpicture}[scale=1.5]
        \tikzset{
          >=stealth',
          pil/.style={
            ->,
            thick,
            shorten <=2pt,
            shorten >=2pt,}
        }

      \begin{scope}
            \draw[very thin,color=gray, opacity=0.3]
              (-0.1,-0.1) grid (4.1,4.1);
            \draw[->, opacity=0.6] (-0.1,0) -- (4.3,0) node[below]
              {$f_1$};
            \draw[->, opacity=0.6] (0,-0.1) -- (0,4.3) node[above]
              {$f_2$};

            \draw[opacity=0.6, dashed, very thick] (0.3,4) -- (4.0,0.5);
            \path (3.5,2.5) node[left, opacity=0.6] {CHIM};

            \draw[very thick, color=black]
              (0.3,4) .. controls (0.9,0.7) .. (4.0,0.5);

            \filldraw[black] (0.3,4) circle (2pt) node[right] {$F_1^*$};
            \filldraw[black] (4.0,0.5) circle (2pt) node[right] {$F_2^*$};

            \filldraw[black, opacity=0.6] (1.04,3.3) circle (2pt);
            \filldraw[black, opacity=0.6] (1.78,2.6) circle (2pt);
            \filldraw[black, opacity=0.6] (2.52,1.9) circle (2pt);
            \filldraw[black, opacity=0.6] (3.26,1.2) circle (2pt);

            \filldraw[black] (0.52,2.8) circle (2pt);
            \filldraw[black] (0.8,1.7) circle (2pt);
            \filldraw[black] (1.5,0.9) circle (2pt);
            \filldraw[black] (2.7,0.6) circle (2pt);

            \path[pil, thick, opacity=0.8]
              (1.04,3.3) edge (0.52,2.8);
            \path[pil, thick, opacity=0.8]
              (1.78,2.6) edge (0.8,1.7);
            \path[pil, thick, opacity=0.8]
              (2.52,1.9) edge (1.5,0.9);
            \path[pil, thick, opacity=0.8]
              (3.26,1.2) edge (2.7,0.6);
			
			\draw[opacity=0.6, dashed, very thick] (0.8,1.7) -- (1.5,0.9) ;
			\draw[opacity=0.6, dashed, very thick] (2.7,0.6) -- (1.5,0.9) ;
			
			\path (1.2,1.4) node[right] {$d_1$};
			\path (2.2,0.95) node[right] {$d_2$};
			
      \end{scope}

  \end{tikzpicture}

%% file: Motivation.tex
Consider a standard bi-objective minimization problem.
\begin{equation}
\label{2OBJ}
\min_{\textbf{x} \in X} \textbf{F}(\textbf{x}) = [F_1(\textbf{x}),
F_2(\textbf{x})]^T
\end{equation}
A straightforward method for constructing a well-sampled discrete Pareto front
representation, as discussed in \cite{Pereyra2009} and \cite{Pereyra2011a},
involves scalarizing the objectives and adding explicit solution spacing
terms to the equality constraint set.

Using a simple weighted sum scalarization, we obtain the following scalar
objective function,
\begin{equation}
\label{2OBJScal}
f(\textbf{x},\lambda) = (1 - \lambda) F_1 (\textbf{x}) + \lambda F_2(\textbf{x})
\mbox{ 
 with }
\begin{array}{rcl} 
 0 \le \lambda \le 1
\end{array}
\end{equation}
the minimizer of which, $\min_{\textbf{x} \in X} f(x, \lambda) = \textbf{f}(\lambda)$, is
a parametric expression of the Pareto front.  The endpoints, $\textbf{f}(0)$ and
$\textbf{f}(1)$ correspond to the minima of the individual objective functions,
$F_1(\textbf{x})$ and $F_2(\textbf{x})$, respectively, while the parameter
$\lambda$ controls the transition between these extremes.

The mapping between $\lambda$ and position along
the Pareto front is often highly nonlinear.  However, as presented in
\cite{Pereyra2011a}, we can exploit this additional degree of freedom to sample
specific points along the curve by treating the $\lambda$ parameter as a design
variable.  This added flexibility is complemented by an auxilary equality
constraint, typically of the form
\begin{equation}
\label{PereyraConst}
||\textbf{F}(\textbf{x}) - \textbf{F}_{prev} ||^2 = \gamma^2
\end{equation}
that explicitly enforces equal spacing between adjacent sample points. This
allows a suitable scalar optimization procedure to minimize the scalarized
objective function \eqref{2OBJScal} subject to the constraint that the result
lies a certain distance, $\gamma$, from a specified point, $\textbf{F}_{prev}$.  By
seeding this procedure with either of the objective minima, for instance with
$\textbf{F}(\textbf{x}_0) = \textbf{f}(0)$, we can iteratively
generate evenly spaced samples by ``marching'' along the Pareto front carrying out the following scalarized minimization
\begin{eqnarray}
\label{PereyraMarching}
\min  \limits_{\textbf{x}_p \in X} & & (1 - \lambda_p) F_1 (\textbf{x}_p) + \lambda_p F_2(\textbf{x}_p) \text{ with }  0 \le \lambda_n \le 1 \\
\text{subject to} & & \nonumber \\
\label{PereyraMarchingConst}
& & ||\textbf{F}(\textbf{x}_p) - \textbf{F}(\textbf{x}_{p-1}) ||^2 = \gamma^2
\end{eqnarray}
for $p$ up to the number of desired sample points, $P$ (provided a suitable
step-size, $\gamma$).

While this method produces well-sampled Pareto front representations, even for
difficult problems, it requires sequential variation of the $\lambda$-parameter
and serial computation of the sample point locations, thereby limiting its
scalability and application to large problems and those with expensive objective
functions.

The authors address this concern in \cite{Pereyra2011a} by extending the
algorithm to parallel processing environments.  This is done, primarily, by
replacing references to previously computed sample points,
$\textbf{F}(\textbf{x}_{p-1})$, with continually updated position estimates,
$\textbf{F}(\textbf{x}_{p-1}^{current})$.  This way, we create a set of
semi-independent and simultaneously solvable scalar optimization sub-problems,
coupled only through the auxilary equispacing constraint \eqref{PereyraConst}.
In fact, we can distribute the subproblems to as many processors as we have
Pareto front sample points, thereby lifting the scalability limits imposed by
the sequential algorithm.  This iterative refinement, of course, requires
defining a set of initial estimates, $\{\textbf{F}_p \for p \in [1,P]\}$,  that
form an ansatz front; however, this can be as simple as a set of equally spaced
points on the line connecting the objective space images of the individual
function minimizers.


One concern with enforcing an explicit spacing constraint like
\eqref{PereyraMarchingConst} is that, without additional information about the arc-length of the true Pareto front, a suitable step-size
parameter, $\gamma$, can not be established a priori.  One can imagine a bootstrapping process whereby the total
arc-length of the front is estimated and refined by successive applications of
the algorithm; however, this would involve serial repetitions of the front
sampling procedure, limiting parallel scalability.



A simpler approach that accomplishes the same goal is to adjust
constraint \eqref{PereyraMarchingConst} to, instead, explicitly enforce
\textit{equal} spacing between a point's up- and down-stream neighbors. 
Specifically, a modified constraint of the form
\begin{equation}
\label{newconst}
||\textbf{F}_p - \textbf{F}_{p+1} ||^2  - ||\textbf{F}_p - \textbf{F}_{p-1}
||^2= 0
\end{equation}
eliminates the need to estimate the Pareto front arc-length a priori by not
fixing the spacing between points and allowing the front representation to
expand indefinitely.  

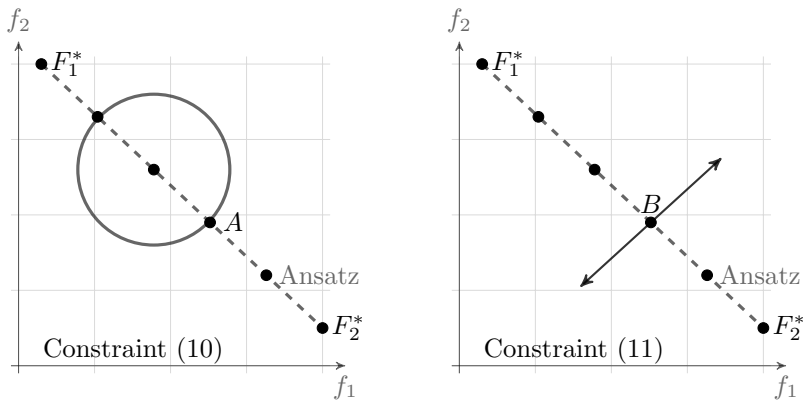
\begin{figure}[h!]
  \centering
  \input{./constraint.tex}
  \caption{The solid circle in the plot on the left represents the feasible
  region permitted by constraint \eqref{PereyraMarchingConst} for ansatz node
  $A$, given unity $\alpha$, while the solid line on the right represents the
  feasible region implied by constraint \eqref{newconst} for ansatz node $B$.}
  \label{fig:constraint}
\end{figure}

Rephrasing the constraint this way, however, has other important ramifications,
as well.  First, it ameliorates an unfortunate consequence of
the formulation in \eqref{PereyraMarchingConst} that limits the optimizer's
improvement per iteration to $\gamma$.  Because of this limited step size, if
the ansatz is sufficiently distant from the real Pareto optimal front, applying
the method of \cite{Pereyra2011a} will require multiple iterations to generate a
set of truly optimal solutions.  Constraint \eqref{newconst}, however, is more in
line with NBI and NC methods (described
above) that do not restrict step sizes, allowing the sample points to converge
to the Pareto front more quickly.  

 
Most importantly, by eliminating explicit references to the total arc-length,
only local communication of current $\textbf{F}_p$ values between adjacent
optimization sub-problems is required.  This means we can concurrently compute
uniformly distributed Pareto front representations using only local,
asynchronous communication between neighboring compute nodes, a key feature of
massively scalable algorithms.  In fact, the amount of parallelism is bounded
only by the number of Pareto front sample points required, $P$, which (to obtain
a truly uniformly spaced representation) often scales exponentially with the
number of objectives.

%% file: constraint.tex
  \begin{tikzpicture}
        \tikzset{
          >=stealth',
          pil/.style={
            ->,
            thick,
            shorten <=2pt,
            shorten >=2pt,}
        }

      \begin{scope}[xshift=-5.8cm]
            \draw[very thin,color=gray, opacity=0.3]
              (-0.1,-0.1) grid (4.1,4.1);
            \draw[->, opacity=0.6] (-0.1,0) -- (4.3,0) node[below]
              {$f_1$};
            \draw[->, opacity=0.6] (0,-0.1) -- (0,4.3) node[above]
              {$f_2$};

            \draw[opacity=0.6, dashed, very thick] (0.3,4) -- (4.0,0.5);
            \path (3.3,1.2) node[right, opacity=0.6] {Ansatz};

            \filldraw[black] (0.3,4) circle (2pt) node[right] {$F_1^*$};
            \filldraw[black] (4.0,0.5) circle (2pt) node[right] {$F_2^*$};

            \filldraw[black] (1.04,3.3) circle (2pt);
            \filldraw[black] (1.78,2.6) circle (2pt);
            \path (2.57,1.9) node[right] {$A$};
            \filldraw[black] (2.52,1.9) circle (2pt);
            \filldraw[black] (3.26,1.2) circle (2pt);

            \draw[very thick, color=black, opacity=0.6] 
				(1.78,2.6) circle (2/2) ;
				
			\path (0.2,0.2) node[right] {Constraint \eqref{PereyraMarchingConst}};

      \end{scope}
      
      \begin{scope}
            \draw[very thin,color=gray, opacity=0.3]
              (-0.1,-0.1) grid (4.1,4.1);
            \draw[->, opacity=0.6] (-0.1,0) -- (4.3,0) node[below]
              {$f_1$};
            \draw[->, opacity=0.6] (0,-0.1) -- (0,4.3) node[above]
              {$f_2$};

            \draw[opacity=0.6, dashed, very thick] (0.3,4) -- (4.0,0.5);
            \path (3.3,1.2) node[right, opacity=0.6] {Ansatz};

            \filldraw[black] (0.3,4) circle (2pt) node[right] {$F_1^*$};
            \filldraw[black] (4.0,0.5) circle (2pt) node[right] {$F_2^*$};

            \filldraw[black] (1.04,3.3) circle (2pt);
            \filldraw[black] (1.78,2.6) circle (2pt);
            \path (2.52,1.9) node[above] {$B$};
            \filldraw[black] (2.52,1.9) circle (2pt);
            \filldraw[black] (3.26,1.2) circle (2pt);

            \path[pil, thick, opacity=0.8]
              (2.52,1.9) edge (1.54,1.0);
            \path[pil, thick, opacity=0.8]
              (2.52,1.9) edge (3.5,2.8);
			
			\path (0.2,0.2) node[right] {Constraint \eqref{newconst}};
			
      \end{scope}

  \end{tikzpicture}

%% file: Generalization.tex
Generalizing this approach to multiple objectives requires not
only auxiliary scalarization variables, but a full generalization of the ansatz
concept and corresponding constraints as well.  The equidistance constraint
naturally applies in the context of the bi-objective linear ansatz.  In fact,
during the first iteration, it reduces the attainable set in a manner that is
functionally equivalent to the constraints used in the NBI and NC methods. 
Similarly confining the attainable set as the dimension of the problem increases
is accomplished by properly balancing the positions of neighboring ansatz nodes
in $k$-space.  In addition to establishing this neighbor-wise
balance, we identified a number of further requirements for any meshing
solution:
\begin{enumerate}[I]
\item \label{MR1} An ansatz mesh must have an approximately equal distribution
of points over its $k-1$ area.  
\item \label{MR2} It should easily support a range of sampling
resolutions.
\item \label{MR3} It should scale with problem dimension. 
\item \label{MR4} 
It should support distributed generation with low computational, as well as
storage, overhead.
\end{enumerate}

To address these requirements, we supply an ansatz with the simplicity and
scalability of a rectilinear, coordinate-parallel mesh that resides on the
convex hull of the individual objective minimizers.  Since we deal in a
normalized objective space, we can generate the initial front sample points by
intersecting a rectilinear grid with the unit $k-1$ simplex.  We then establish
equidistance constraint relationships (or adjacency pairs) along the
coordinate-parallel mesh axes.  For instance, the scalarized sub-problem at
point $\textbf{C}$ will be augmented with two additional constraints,
\eqref{Cvertical} and \eqref{Chorizontal}, both of which are based on
\eqref{newconst}.  As in the two objective problem discussed above, this
initially reduces the attainable set of this optimization subproblem to the
simplex normal.

\begin{equation}
\label{Cvertical}
||\textbf{F}_C - \textbf{F}_{A} ||^2  - ||\textbf{F}_C - \textbf{F}_{E}
||^2= 0
\end{equation}
\begin{equation}
\label{Chorizontal}
||\textbf{F}_C - \textbf{F}_{B} ||^2  - ||\textbf{F}_C - \textbf{F}_{D}
||^2= 0
\end{equation}

The scalarized subproblem at point $\textbf{D}$, however, will only be augmented with a single additional constraint, \eqref{Ddiag}.  This is because point $\textbf{D}$ actually lies on the Pareto front of the bi-objective subproblem in the $f_1 \times f_2$ plane and needs no further restriction of the attainable set.

\begin{equation}
\label{Ddiag}
||\textbf{F}_D - \textbf{F}_{A} ||^2  - ||\textbf{F}_D - \textbf{F}_1^*
||^2= 0
\end{equation}

\begin{figure}[h!]
  \centering
  \input{./3dAnsatz.tex}
  
  \caption{The full ansatz for a three objective problem as seen in the $f_1
  \times f_2 \times f_3$ space.}
  \label{fig:3dansatz}
\end{figure}
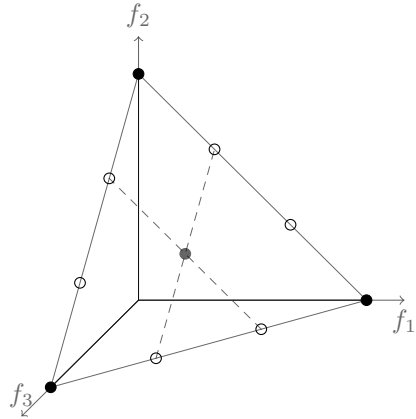

The ansatz is completed by reprojecting the
result into $k$-space resulting in a uniformly-distributed linear interpolation
between the $k$ individual objective minimizers and preserving the
coordinate-wise equal spacing between the constituent points. This process is
illustrated for a three objective problem in figure \ref{fig:ansatzConstruction}
with the final result projected back into 3-space in figure \ref{fig:3dansatz}.

\begin{figure}[h!]
  \centering
  \input{./ansatzConstruction.tex}
  
  \caption{On the left, we generate the ansatz points by intersecting a  
  rectilinear grid with the convex hull of the individual minimizers, denoted by
  the shaded region (a $k-1$ simplex).  We establish neighbor-relations on the
  right using coordinate-parallel axes.}
  \label{fig:ansatzConstruction}
\end{figure}
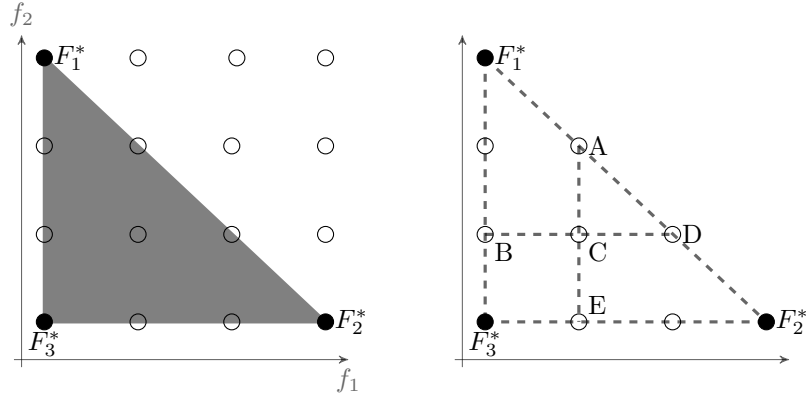


This ansatz concept addresses requirements \ref{MR1} and \ref{MR2} by creating a
uniform tensor-product mesh and supporting a range of axis-wise sampling densities.  In
this way, we avoid uneven sampling of the Pareto front and allow the user to
explicitly specify the level of detail required.
 
Moreover, the coordinate-parallel nature of this solution addresses
characteristic \ref{MR3} by adding only two more mesh point neighbors per additional dimension
of the ansatz.  In other words, the degree of each node in the graph is bounded
by $2(k-1)$ while adequately reducing the attainable set.  Ultimately, this means
that the number of additional constraints added to the overarching optimization problem to ensure equispacing of the Pareto front, as
well as the magnitude of the data dependency of the optimization sub-problem
occurring at each mesh node, scales linearly with the number of target objectives.

%

\begin{figure}[h!]
  \centering
  \input{./4dAnsatz.tex}
  
  \caption{A similarly constructed ansatz for a four objective problem,
  projected onto the $f_1 \times f_2 \times f_3$ subspace.}
  \label{fig:4dansatz}
\end{figure}
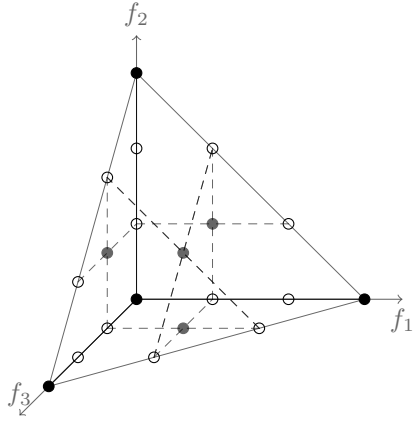

Finally, since all of the coordinate transforms and neighbor-discovering
operations can be accomplished independently for each point, this meshing
strategy supports fully distributed generation, satisfying requirement
\ref{MR4}.

%

With a scalable ansatz concept in place, the fully generalized approach boils
down to solving a scalar minimization problem, given by \eqref{alg:Homotopy},
for each Pareto front sample point, $p$.  Here, an adjacency pair, $a_p$,
consists of a set of per-axis opposing neighbors of point $p$ with $a_p^l$ and
$a_p^r$ denoting the mesh indicies of the left and right neigbors, respectively.
$A(p)$ represents the set of all adjacency pairs of ansatz mesh point $p$.

\begin{eqnarray}
\label{alg:Homotopy}
\min\limits_{\textbf{x} \in X, \boldsymbol\lambda} \lambda_1 F_1 (\textbf{x}) &+& \lambda_2 F_2(\textbf{x}) + \ldots + \lambda_k F_k(\textbf{x}) \text{ with } 0 \le \lambda_i \le 1 \\
\text{ subject to } & & \nonumber \\
g_j(\textbf{x}) &\le & 0 \quad \forall j \in [1, m] \nonumber \\
h_l(\textbf{x}) &=& 0 \quad \forall l \in [1, e] \nonumber \\
\sum_{i=1}^k\lambda_i &=& 1 \mbox{, } \quad  \forall i \in [1,k] \nonumber \\
||\textbf{F}(\textbf{x}) &-& \textbf{F}_{a_p^l} ||^2  - ||\textbf{F}(\textbf{x}) - \textbf{F}_{a_p^r} ||^2 = 0 \quad  \forall a_p \in A(p) \nonumber
\end{eqnarray}

%
%

%% file: ansatzConstruction.tex
  \begin{tikzpicture}
        \tikzset{
          >=stealth',
          pil/.style={
            ->,
            thick,
            shorten <=2pt,
            shorten >=2pt,}
        }

        \begin{scope}[xshift=-5.8cm]
            \draw[->, opacity=0.6] (-0.1,0) -- (4.3,0) node[below]
              {$f_1$};
            \draw[->, opacity=0.6] (0,-0.1) -- (0,4.3) node[above]
              {$f_2$};
              
            \filldraw[very thick, color=black,opacity=0.5] (0.3,4) -- (0.3,0.5)
            -- (4.0,0.5) -- cycle; 

            \filldraw[black] (0.3,4) circle (3pt) node[right] {$F_1^*$};
            
            \draw[black] (1.533,2.833) circle (3pt);
            \draw[black] (2.766,1.66) circle (3pt);
            
            \filldraw[black] (4.0,0.5) circle (3pt) node[right] {$F_2^*$};
            
            \draw[black] (0.3,2.833) circle (3pt);
            \draw[black] (0.3,1.66) circle (3pt);
            
            \filldraw[black] (0.3,0.5) circle (3pt) node[below] {$F_3^*$};
            
            \draw[black] (1.533,1.66) circle (3pt);
            \draw[black] (1.533,0.5) circle (3pt);
            
            \draw[black] (2.766,0.5) circle (3pt);
			
			\draw[black] (1.533,4) circle (3pt);
			\draw[black] (2.833,4) circle (3pt);
			\draw[black] (4,4) circle (3pt);
            
            \draw[black] (2.766,2.833) circle (3pt);
            \draw[black] (4,2.833) circle (3pt);
            \draw[black] (4,1.66) circle (3pt);

			
      \end{scope}

      \begin{scope}
            \draw[->, opacity=0.6] (-0.1,0) -- (4.3,0) ;
            \draw[->, opacity=0.6] (0,-0.1) -- (0,4.3) ;
              
            \draw[very thick, color=black, opacity=0.6, dashed] 
				(0.3,4) -- (0.3,0.5);
			\draw[very thick, color=black, opacity=0.6, dashed] 
				(0.3,4) -- (4.0,0.5);
			\draw[very thick, color=black, opacity=0.6, dashed] 
				(0.3,0.5) -- (4.0,0.5);
			
			\draw[very thick, color=black, opacity=0.6, dashed] 
				(0.3,1.66) -- (2.766,1.66);
			\draw[very thick, color=black, opacity=0.6, dashed] 
				(1.533,2.833) -- (1.533,0.5);

            \filldraw[black] (0.3,4) circle (3pt) node[right] {$F_1^*$};
            
            \draw[black] (1.533,2.833) circle (3pt) node[right] {A};
            \draw[black] (2.766,1.66) circle (3pt) node[right] {D};
            
            \filldraw[black] (4.0,0.5) circle (3pt) node[right] {$F_2^*$};
            
            \draw[black] (0.3,2.833) circle (3pt);
            \draw[black] (0.3,1.66) circle (3pt) ;
            \draw[black] (0.3,1.66-.23) node[right] {B};
            
            \filldraw[black] (0.3,0.5) circle (3pt) node[below] {$F_3^*$};
            
            \draw[black] (1.533,1.66) circle (3pt);
            \draw[black] (1.533,1.66-.23) node[right] {C};
            \draw[black] (1.533,0.5) circle (3pt);
            \draw[black] (1.533,0.5 +.23)  node[right] {E};
            
            \draw[black] (2.766,0.5) circle (3pt);

      \end{scope}
  \end{tikzpicture}

%% file: 3dAnsatz.tex
\begin{tikzpicture}[line join=round]%

\draw[->, opacity=0.6] (0,0) -- (3.5,0) node[below]
              {$f_1$};
\draw[->, opacity=0.6] (0,0) -- (0,3.5) node[above]
              {$f_2$};
\draw[->, opacity=0.6] (0,0,0) -- (0,0,4) node[above]
              {$f_3$};

    \coordinate (A) at ( 0,       0,     0);%
    \coordinate (B) at (1.0,    0,     0);%
    \coordinate (C) at (2*1.0,  0,     0);%
    \coordinate (D) at (3*1.0,  0,     0);%
    
    \coordinate (E) at ( 0,       1.0,     0);%
    \coordinate (F) at (1.0,    1.0,     0);%
    \coordinate (G) at (2*1.0,  1.0,     0);%
    
    \coordinate (H) at ( 0,       2*1.0,     0);%
    \coordinate (I) at (1.0,    2*1.0,     0);%
    
    \coordinate (J) at (0,    3*1.0,     0);%
    \coordinate (AA) at ( 0,       0,     1.0);%
    \coordinate (BB) at (1.0,    0,     1.0);%
    \coordinate (CC) at (2*1.0,  0,     1.0);%
    
    \coordinate (DD) at ( 0,       1.0,     1.0);%
    \coordinate (EE) at (1.0,    1.0,     1.0);%
    
    \coordinate (FF) at (0,  2*1.0,     1.0);%
    
    \coordinate (AAA) at ( 0,       0,     2*1.0);%
    \coordinate (BBB) at (1.0,    0,     2*1.0);%
    
    \coordinate (CCC) at (0,    1.0,     2*1.0);%

	\coordinate (AAAA) at (0,    0,     3*1.0);%
    
    \draw [very thin, black,  opacity=0.6, line cap=rect]%
        (A) -- (J) -- (D) -- cycle;%
        
    \draw [very thin, black,  opacity=0.6, line cap=rect]%
        (A) -- (D) -- (AAAA) -- cycle;%
        
    \draw [very thin, black,  opacity=0.6, line cap=rect]%
        (A) -- (J) -- (AAAA) -- cycle;%


    \draw [very thin, black, dashed, opacity=0.6, line cap=rect]%
        (BBB) -- (I);%
    \draw [very thin, black, dashed, opacity=0.6, line cap=rect]%
        (CC) -- (FF);%

    \filldraw[black] (D) circle (2pt);

    \draw[black] (G) circle (2pt);
    \draw[black] (I) circle (2pt);
    \filldraw[black] (J) circle (2pt);

    \draw[black] (CC) circle (2pt);
    \filldraw[black, opacity=0.6] (EE) circle (2pt);
    \draw[black] (FF) circle (2pt);

    \draw[black] (BBB) circle (2pt);
    \draw[black] (CCC) circle (2pt);
    
    \filldraw[black] (AAAA) circle (2pt);
    
\end{tikzpicture}%

%% file: 4dAnsatz.tex
\begin{tikzpicture}[line join=round]%

\draw[->, opacity=0.6] (0,0) -- (3.5,0) node[below]
              {$f_1$};
\draw[->, opacity=0.6] (0,0) -- (0,3.5) node[above]
              {$f_2$};
\draw[->, opacity=0.6] (0,0,0) -- (0,0,4) node[above]
              {$f_3$};

    \coordinate (A) at ( 0,       0,     0);%
    \coordinate (B) at (1.0,    0,     0);%
    \coordinate (C) at (2*1.0,  0,     0);%
    \coordinate (D) at (3*1.0,  0,     0);%
    
    \coordinate (E) at ( 0,       1.0,     0);%
    \coordinate (F) at (1.0,    1.0,     0);%
    \coordinate (G) at (2*1.0,  1.0,     0);%
    
    \coordinate (H) at ( 0,       2*1.0,     0);%
    \coordinate (I) at (1.0,    2*1.0,     0);%
    
    \coordinate (J) at (0,    3*1.0,     0);%
    \coordinate (AA) at ( 0,       0,     1.0);%
    \coordinate (BB) at (1.0,    0,     1.0);%
    \coordinate (CC) at (2*1.0,  0,     1.0);%
    
    \coordinate (DD) at ( 0,       1.0,     1.0);%
    \coordinate (EE) at (1.0,    1.0,     1.0);%
    
    \coordinate (FF) at (0,  2*1.0,     1.0);%
    
    \coordinate (AAA) at ( 0,       0,     2*1.0);%
    \coordinate (BBB) at (1.0,    0,     2*1.0);%
    
    \coordinate (CCC) at (0,    1.0,     2*1.0);%

	\coordinate (AAAA) at (0,    0,     3*1.0);%
    
    \draw [very thin, black,  opacity=0.6, line cap=rect]%
        (A) -- (J) -- (D) -- cycle;%
        
    \draw [very thin, black,  opacity=0.6, line cap=rect]%
        (A) -- (D) -- (AAAA) -- cycle;%
        
    \draw [very thin, black,  opacity=0.6, line cap=rect]%
        (A) -- (J) -- (AAAA) -- cycle;%

    \draw [very thin, black, dashed, opacity=0.6, line cap=rect]%
        (B) -- (I);%
    \draw [very thin, black, dashed, opacity=0.6, line cap=rect]%
        (E) -- (G);%
        
    \draw [very thin, black, dashed, opacity=0.6, line cap=rect]%
        (BBB) -- (I);%
    \draw [very thin, black, dashed, opacity=0.6, line cap=rect]%
        (CC) -- (FF);%
        
    \draw [very thin, black, dashed, opacity=0.6, line cap=rect]%
        (BBB) -- (I);%
    \draw [very thin, black, dashed, opacity=0.6, line cap=rect]%
        (CC) -- (FF);%
        
    \draw [very thin, black, dashed, opacity=0.6, line cap=rect]%
        (AA) -- (FF);%
    \draw [very thin, black, dashed, opacity=0.6, line cap=rect]%
        (CCC) -- (E);%
        
    \draw [very thin, black, dashed, opacity=0.6, line cap=rect]%
        (BBB) -- (B);%
    \draw [very thin, black, dashed, opacity=0.6, line cap=rect]%
        (AA) -- (CC);%
    
    \filldraw[black] (A) circle (2pt);
    \draw[black] (B) circle (2pt);
    \draw[black] (C) circle (2pt);
    \filldraw[black] (D) circle (2pt);
    \draw[black] (E) circle (2pt);
    \filldraw[black, opacity=0.6] (F) circle (2pt);
    \draw[black] (G) circle (2pt);
    \draw[black] (H) circle (2pt);
    \draw[black] (I) circle (2pt);
    \filldraw[black] (J) circle (2pt);
    
    \draw[black] (AA) circle (2pt);
    \filldraw[black, opacity=0.6] (BB) circle (2pt);
    \draw[black] (CC) circle (2pt);
    \filldraw[black, opacity=0.6] (DD) circle (2pt);
    \filldraw[black, opacity=0.6] (EE) circle (2pt);
    \draw[black] (FF) circle (2pt);

	\draw[black] (AAA) circle (2pt);
    \draw[black] (BBB) circle (2pt);
    \draw[black] (CCC) circle (2pt);
    
    \filldraw[black] (AAAA) circle (2pt);
    
\end{tikzpicture}%

%% file: Results.tex
\subsection{Metrics}

Having described the new algorithm (\ref{alg:Homotopy}), we
turn to evaluating its performance relative to the other methods described.  We discuss a set of metrics for measuring front accuracy as well as quality (with a focus on evenly sampled representations) and compute them for a set of multi-objective optimization benchmark problems taken from
the literature.

\subsection{Hypervolume}

A useful and, more importantly, scalable metric for measuring the accuracy and
quality of a Pareto front approximation is the hypervolume of the dominated
space. Effectively, given a point in the attainable set, we measure the $k$
dimensional volume of the dominated space, relative to some ``worst-case''
point, as shown in figure \ref{fig:metricsDiagram}.
 
The particular convenience of this metric lies in the fact that the limit of a
perfectly sampled Pareto front will converge to a concrete value, representing
the total dominated volume.  Given the fact that the true Pareto optimal front
is the maximal set of the attainable region poset, with this metric, we are
able to directly measure convergence to the true Pareto front.  We employed the
implementation provided by \cite{While2012} for efficient computation of
hypervolumes.  

\subsection{Evenness}

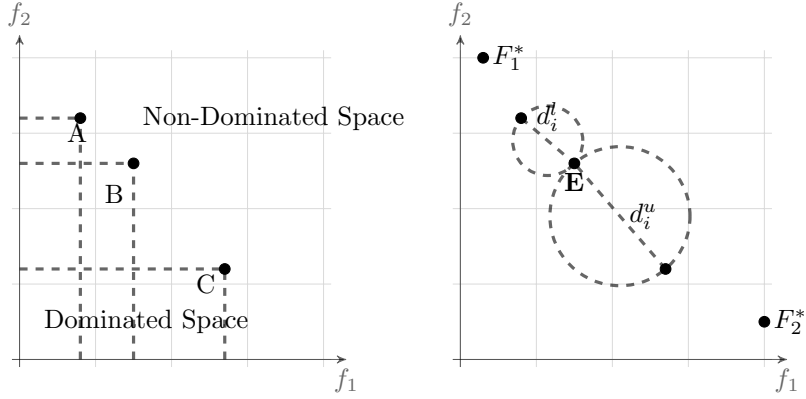
\begin{figure}[h!]
  \centering
  \input{./metrics.tex}
  
  \caption{The hypervolume of a Pareto front is the total volume of the
  dominated space, as illustrated on the left.  Computed with respect to the
  origin, this equals the volume of the union of sets A, B, and C.  With regard to the evenness metric, the plot on the right demonstrates the construction of $d_i^l$ and $d_i^u$ for a point $\textbf{E}$.}
  \label{fig:metricsDiagram}
\end{figure}

Beyond simply converging to the Pareto optimal front, we aim to ensure a
uniformly sampled approximation in an effort to maximize both economy of
function evaluations as well as permit useful interpretation of results.  A
number of methods have been proposed to capture the essence of uniformity, from
normalized root mean square (RMS) distances between adjacent points
\cite{Pereyra2009} to the maximum lower bound of all inter-point distances \cite{Leyffer2008}; however,
these methods do not generalize to varying objective scales and are not well
adopted by the community at large.  In \cite{Messac2004}, the authors recommend
a measure of ``evenness'' that corresponds to the intuitive notion that no region of
the Pareto front is over- or under-represented by the discrete approximation. 
Given a point $F_i$ in the approximation set, we first consider the nearest
neighbor distance, $d_i^l$, to another point in the set, as in figure
\ref{fig:metricsDiagram}.  Next, we consider the largest sphere, containing no points
from the set, that can be constructed such that $F_i$ and another point, $F_j$,
both lie on the surface and denote the diameter, $d_i^u$.  Then, given the set
\begin{equation}
d \colon= \left\{ d_i^l, d_i^u \colon \forall i \in [1, S] \right\}
\end{equation}
we define the evenness, $\textbf{E}$, as  
\begin{equation}
\textbf{E} = \sigma_d / \hat{d}
\end{equation}
the ratio of this set's standard deviation, $\sigma_d$, to its mean, $\hat{d}$,
implying that a perfectly uniform distribtion has $\textbf{E} = 0$.  This metric
generalizes to multiple dimensions quite easily by considering
hyperspheres when determining the $d^l$ and $d^u$ values. 

\subsection{Function Evaluations \& Time-to-solution}

Finally, while the previous two metrics quantify front accuracy and quality, we
aim to measure the amount of computational effort required to obtain the given
solutions.  To directly compare algorithms, independent of specific
implementation or computing platform, we report the number of objective function
evaluations per Pareto point used to generate front approximations for the
synthetic test problems.  


\subsection{Test Problems}

Since the homotopy method presented here primarily relies on augmenting the
constraint functions, $\textbf{h}(\textbf{x})$, with additional equispacing
constraints, we chose to directly compare the accuracy and quality of the Pareto
front approximation with those produced by similar methods, like NBI and NC. 
Therefore, we selected a set of test problems taken from the recent
literature \cite{S.Motta2012} and where some of the above metrics are reported
for a variety of scalarization and geometric methods.

\subsection{Problem 1}

The first problem is relatively easy, consisting of linear functions and convex
constraints.  
\begin{eqnarray}
\label{Motta1}
\min\limits_{\textbf{x}} f_i(\textbf{x}) &=& x_i \mbox{ for } i = 1, 2, 3 \\
\text{ subject to } & & \nonumber \\
& & -x_1 + x_2^{-1} + x_3^{-1} \le 0 \nonumber \\
& & x_1^{-1} - x_2 + x_3^{-1} \le 0 \nonumber \\
& & x_1^{-1} + x_2^{-1} - x_3 \le 0 \nonumber \\
& & 0.2 \le x_i \le 10 \nonumber
\end{eqnarray}

\begin{figure}[h!]
  \centering
  \includegraphics[scale=0.4]{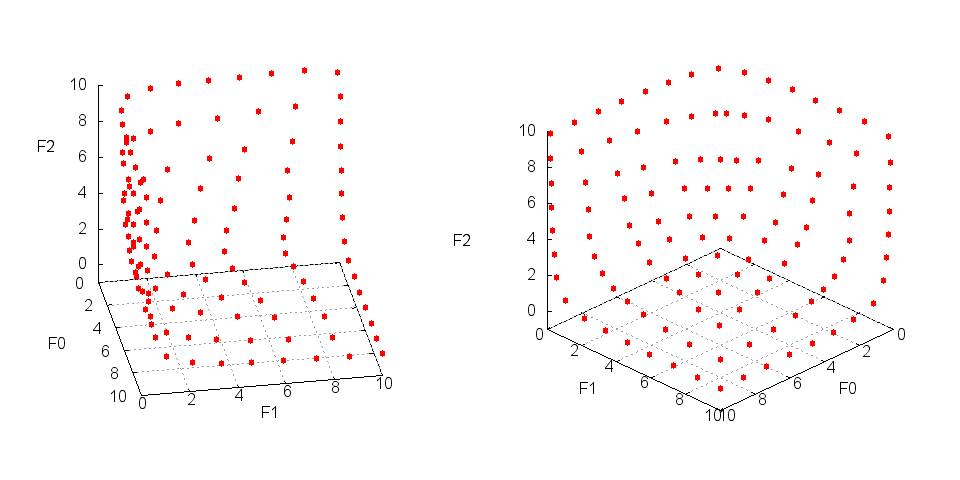}
  \caption{Pareto front for problem \eqref{Motta1} consisting of 120
  points obtained with 6,264 function evaluations.}
  \label{graph:Motta1}
\end{figure}

Using the simplical CHIM 120-point ansatz, we obtained a smooth, connected
Pareto front (figure \ref{graph:Motta1}) covering most of the first octant using four
iterations of our algorithm and a few thousand function evaluations.

We computed evenness values, as well as state the number of
function evaluations required per Pareto point, and compare them to reported
values from the literature in table \ref{tab:Motta1}.

\bigskip
\begin{table}
\caption{Problem \eqref{Motta1} Results, A comparison of solution quality (evenness [\textbf{E}]) and computational
effort (function evaluations per point [\textbf{FE/Point}]) for problem \eqref{Motta1}.
Reference values for WS, NBIm, and NCm algorithms were obtained from
\cite{S.Motta2012}.  It is not clear whether the number of function evaluations
reported includes those induced by finite difference schemes, or if
analytical derivatives were used in the scalar minimizations.  Our results were
obtained with a first-order finite difference implementation.} \label{tab:Motta1} 
\centering
\begin{tabular}{c c c}
\hline \hline 
\textbf{Method} &  $\textbf{E}$ & \textbf{FE/Point} \\ [0.5ex] 
\hline
\\  WS       &  4.331     &  120.2     \\
\\  NBIm     &  0.2958     &  34.3     \\
\\  NCm      &  0.2958     &  34.4     \\
\\  Our new algorithm (\ref{alg:Homotopy}) & 0.01930      &  52.2   \\ [1ex]
\hline
\end{tabular}
\par
\end{table}

\subsection{Problem 2}

The other problem is an extension of \eqref{Motta1} to four objectives.

\begin{eqnarray}
\label{Motta2}
\min\limits_{\textbf{x}} f_i(\textbf{x}) &=& x_i \mbox{ for } i = 1, 2, 3, 4 \\
\text{ subject to } & & \nonumber \\
& & -x_1 + x_2^{-1} + x_3^{-1} + x_4^{-1} \le 0 \nonumber \\
& & x_1^{-1} - x_2 + x_3^{-1} + x_4^{-1}  \le 0 \nonumber \\
& & x_1^{-1} + x_2^{-1} - x_3 + x_4^{-1}  \le 0 \nonumber \\
& & x_1^{-1} + x_2^{-1} + x_3^{-1} - x_4 \le 0 \\
& & 0.2 \le x_i \le 10 \nonumber
\end{eqnarray}

\begin{figure}[h!]
  \centering
  \includegraphics[scale=0.4]{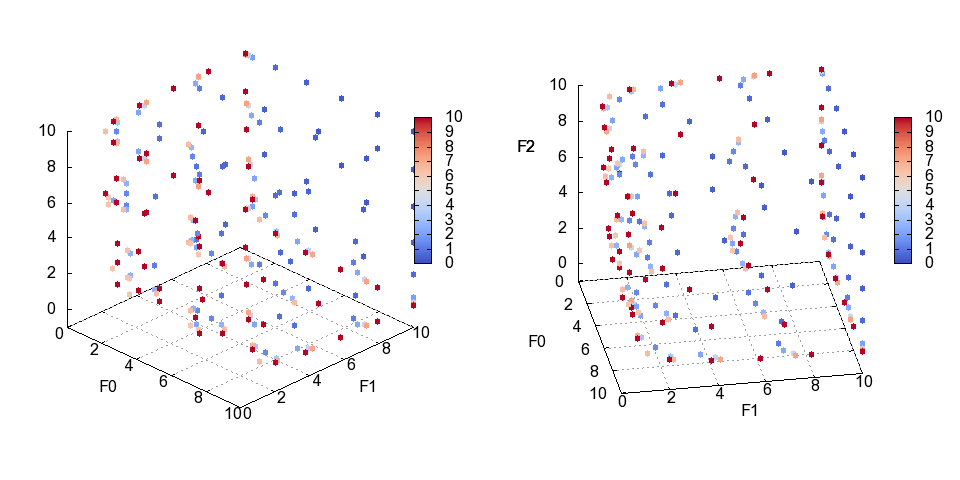}
  \caption{Pareto front for problem \eqref{Motta2} consisting of 220
  points obtained with 10,670 function evaluations.  The value of the fourth
  objective for each point is shown using the red-blue color-map.}
  \label{graph:Motta2}
\end{figure}

Likewise, with a simplicial CHIM 220-point ansatz, we obtained a smooth,
connected Pareto front (figure \ref{graph:Motta2}) using two
iterations of our algorithm.

As before, we computed evenness values, and tracked the
number of function evaluations required per Pareto point.\ Table
\ref{tab:Motta2} compares the results to reported values from the literature.

\bigskip
\begin{table}
\caption{Problem \eqref{Motta2} Results  A comparison of solution quality (evenness [\textbf{E}]) and computational
effort (function evaluations per point [\textbf{FE/Point}]) for problem \eqref{Motta2}.
Reference values for WS, NBIm, and NCm algorithms were obtained from
\cite{S.Motta2012}. } \label{tab:Motta2},
\centering
\begin{tabular}{c c c}
\hline \hline 
\textbf{Method} &  $\textbf{E}$ & \textbf{FE/Point} \\ [0.5ex]  
\hline
   WS	           &  4.836     &  231.9     \\
\\ NBIm           &  0.3262     &  49.3     \\
\\ NCm            &  0.3072     &  55.3     \\
\\ Our new algorithm (\ref{alg:Homotopy})      & 0.02091     &  48.5 \\ [0.5ex]
\hline
\end{tabular}
\par
\end{table}

\subsection{Problem 3}

On the other hand, the EA community tends to focus on different kinds of
problems.  EAs, being driven by random generation and recombination of
solutions, are typically useful in situations where derivative information is
unavailable or unhelpful.  As a result, test problem suites intended for EAs
employ highly non-convex, rapidly oscillating, and extremely discontinuous
objectives that are often inappropriate for constraint-based methods, like the
one presented here.  Therefore, we limit our discussion to a representative
example as the primary motivation is to demonstrate the precise equispacing of
produced solutions and economy of function evaluations.

\begin{eqnarray}x
\label{DTLZ2}
\min\limits_{\textbf{x}} & & f_i(\textbf{x}) \mbox{ for } i = 1, 2, 3 \\ \text{ subject to } \nonumber \\
f_1(\textbf{x}) &= & (1 + g( x ))  \cos( x_1 \pi/2 )  \cos( x_2 \pi/2 ) \nonumber\\
f_2(\textbf{x}) &= & (1 + g( x ))  \cos( x_1 \pi/2 )  \sin( x_2 \pi/2 ) \nonumber \\
f_3(\textbf{x}) &= & (1 + g( x ))  \sin( x_1 \pi/2 ) \nonumber \\
g(\textbf{x}) &= & \sum_{i=2}^d (x_i - 0.5)^2 \mbox{ with } 0 \le x_i \le 1 \mbox{ for } i=1,2 \nonumber
\end{eqnarray}

We selected a workable problem, defined in \eqref{DTLZ2}, from a standard EA test
suite \cite{Deb2005} and compare the performance to some other heuristics noted in the literature.

\begin{figure}[h!]
  \centering
  \includegraphics[scale=0.4]{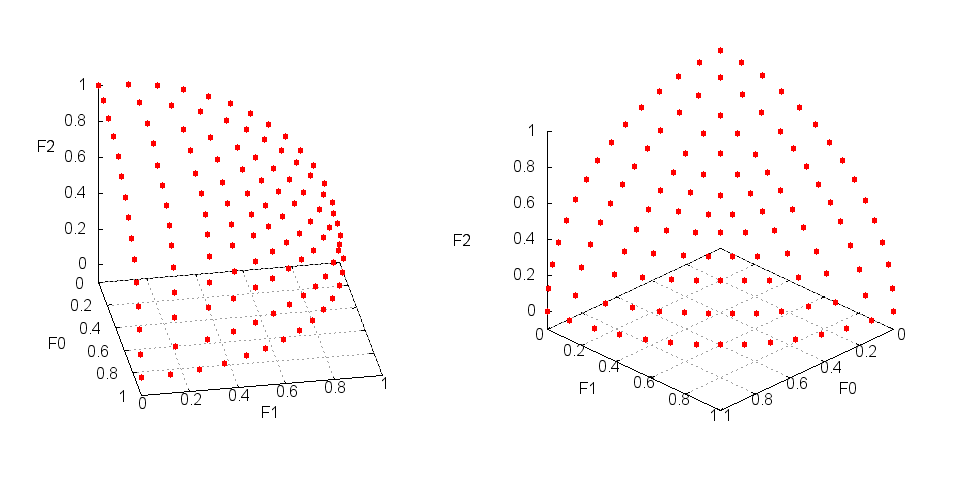}
  \caption{Pareto front for problem \eqref{DTLZ2} consisting of 120
  points obtained with 15,765 function evaluations.  The finite differences
  scheme coupled with the fact that the problem has 10 independent variables
  gives rise to the very high number of function evaluations required.  For
  examples of typical fronts generated by NSGA-II and SPEA2, see
  \cite{Deb2005}.}
  \label{graph:DTLZ2}
\end{figure}

The resulting Pareto front, as seen if figure \ref{graph:DTLZ2}, is computed
using significantly fewer function evalutations, obtains a much better
hypervolume than any of the algorithms compared in either \cite{Nebro2008} or
\cite{Li2012}, and displays more regular sampling than \cite{Deb2005}.

\begin{figure}[h!]
  \centering
  \includegraphics[scale=0.5]{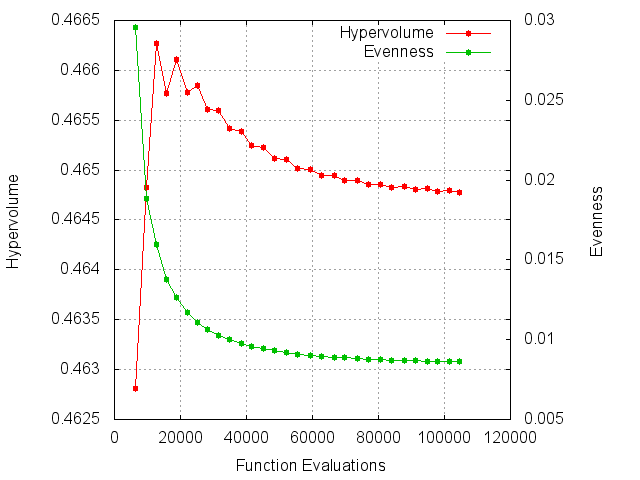}
  \caption{Solution quality metrics for Pareto front for problem
  \eqref{DTLZ2}.  Both the solution hypervolume and evenness improve
  asymptotically as the number of function evaluations increase.  In this
  example, one iteration of algorithm \ref{alg:Homotopy} corresponds to about
  3,000 function evaluations.  The front hypervolume is computed with respect
  to the origin.}
  \label{graph:metrics}
\end{figure}

In figure \ref{graph:metrics}, we demonstrate the asymptotic
improvement in both hypervolume (computed with respect to the origin) and
evenness.  This implies that, for certain simple front topologies, algorithm
\ref{alg:Homotopy} can be terminated after only a few iterations.  Moreover,
figure \ref{graph:metrics} shows that further iterations of the algorithm tend
to improve the spacing of points more than the obtained hypervolume.  This is a
result of directly running a set of scalar optimizations (allowing samples to
converge to the true optimal front more quickly) rather than ranking and
recombining randomly generated solutions.

%% file: metrics.tex
  \begin{tikzpicture}
        \tikzset{
          >=stealth',
          pil/.style={
            ->,
            thick,
            shorten <=2pt,
            shorten >=2pt,}
        }

      \begin{scope}[xshift=-5.8cm]
            \draw[very thin,color=gray, opacity=0.3]
              (-0.1,-0.1) grid (4.1,4.1);
            \draw[->, opacity=0.6] (-0.1,0) -- (4.3,0) node[below]
              {$f_1$};
            \draw[->, opacity=0.6] (0,-0.1) -- (0,4.3) node[above]
              {$f_2$};


            \filldraw[black] (0.8,3.2) circle (2pt);
            \draw[very thick, color=black, opacity=0.6, dashed] 
				(0.0,3.2) -- (0.8,3.2);
			\draw[very thick, color=black, opacity=0.6, dashed] 
				(0.8,0.0) -- (0.8,3.2);
            
            \filldraw[black] (1.5,2.6) circle (2pt);
            \draw[very thick, color=black, opacity=0.6, dashed] 
				(0.0,2.6) -- (1.5,2.6);
			\draw[very thick, color=black, opacity=0.6, dashed] 
				(1.5,0.0) -- (1.5,2.6);

            \filldraw[black] (2.7,1.2) circle (2pt);
			\draw[very thick, color=black, opacity=0.6, dashed] 
				(0.0,1.2) -- (2.7,1.2);
			\draw[very thick, color=black, opacity=0.6, dashed] 
				(2.7,0.0) -- (2.7,1.2);
				
			\path (0.2,0.5) node[right] {Dominated Space};
			\path (1.5,3.2) node[right] {Non-Dominated Space};
			
			\path (0.5,3.0) node[right] {A};
			\path (1.0,2.2) node[right] {B};
			\path (2.2,1.0) node[right] {C};
			
      \end{scope}

      \begin{scope}
            \draw[very thin,color=gray, opacity=0.3]
              (-0.1,-0.1) grid (4.1,4.1);
            \draw[->, opacity=0.6] (-0.1,0) -- (4.3,0) node[below]
              {$f_1$};
            \draw[->, opacity=0.6] (0,-0.1) -- (0,4.3) node[above]
              {$f_2$};

            \filldraw[black] (0.3,4) circle (2pt) node[right] {$F_1^*$};
            \filldraw[black] (4.0,0.5) circle (2pt) node[right] {$F_2^*$};

            \filldraw[black] (0.8,3.2) circle (2pt);
            
            \filldraw[black] (1.5,2.6) circle (2pt) node [below]
            {$\textbf{E}$};
            \draw[very thick, color=black, opacity=0.6, dashed] 
				(1.15,2.9) circle (.922/2);
			\draw[very thick, color=black, opacity=0.6, dashed] 
				(0.8,3.2) -- (1.5,2.6);
			\draw[very thick, color=black, opacity=0.6, dashed] 
				(2.1,1.9) circle (1.844/2) ;
			\draw[very thick, color=black, opacity=0.6, dashed] 
				(2.7,1.2) -- (1.5,2.6);
	
            \filldraw[black] (2.7,1.2) circle (2pt);
				
			\path (2.1,1.9) node[right] {$d_i^u$};
			\path (1.15,2.9) node[above] {$d_i^l$};

      \end{scope}

  \end{tikzpicture}

%% file: Conclusions.tex
Here, an new algorithm (\ref{alg:Homotopy}) for multi-objective optimization
that produces uniformly sampled Pareto fronts by construction is presented. 
While the algorithm is general, it is most suitable for application to
simulation-based engineering optimization problems where economy of function
evaluations and smoothness of result are critical.  The algorithm discussed
achieves an order of magnitude improvement over other geometrically motivated
methods, like Normal Boundary Intersection and Normal Constraint, with respect to solution
evenness for similar computational expense.  This benefit of the proposed method
remains, and even improves, after scaling the number of dimensions (and
therefore the difficulty of the problem).  Furthermore, the resulting
uniformity of solutions extends even to more difficult problems, such as those
appearing in common EA test cases.  

While the resulting discrete representation
of the Pareto front, and computational expense of achieving it, are both improved, an other important
aspect of the proposed method is its amenability to parallelization. We will report on the parallel aspects together with real world problems \cite{ineichen2011} in a forthcoming paper.

%% file: Bibliography.tex
\bibliographystyle{plain}
\bibliography{Thesis}